\documentclass[11pt]{article}
\usepackage{mathrsfs}
\usepackage{amssymb}
\usepackage{latexsym}
\usepackage{amsmath,amsthm}
\usepackage{amsfonts}
\usepackage{url}

\newcommand{\h}{\eta}

\usepackage{bbm}
\newcommand{\chr}{\boldsymbol{\mathbbm{1}}} 
\newcommand{\pred}[1]{\chr_{\left\{ #1 \right\}}}
\newcommand{\mexp}{{\mathbb E}}
\newcommand{\TV}[1]{\nrm{#1}_{\textrm{{\tiny \textup{TV}}}}}
\newcommand{\Lip}[1]{\nrm{#1}_{\textrm{{\tiny \textup{Lip}}}}}
\renewcommand{\P}{\prs}
\newcommand{\prs}{{\mathbb P}}
\newcommand{\pr}[1]{\prs\!\tlprn{#1}}

\newcommand{\inex}{\alpha}

\newcommand{\one}{\text{i}}
\newcommand{\two}{\text{ii}}

\newcommand{\iia}{\ensuremath{\mathrm{(a)}}}
\newcommand{\iib}{\ensuremath{\mathrm{(b)}}}
\newcommand{\iic}{\ensuremath{\mathrm{(c)}}}
\newcommand{\iid}{\ensuremath{\mathrm{(d)}}}
\newcommand{\iie}{\ensuremath{\mathrm{(e)}}}
\newcommand{\iif}{\ensuremath{\mathrm{(f)}}}

\newcommand{\sgn}{\operatorname{sgn}}

\usepackage{geometry}
\makeatletter
\theoremstyle{plain}
\newtheorem{thm}{Theorem}[section]
 \theoremstyle{remark}
 \newtheorem{rem}[thm]{Remark}
 \theoremstyle{plain}    
 \newtheorem{lem}[thm]{Lemma} 
 \theoremstyle{plain}    
 \newtheorem{cor}[thm]{Corollary} 
\makeatother
\newcommand{\bethn}{\begin{thm}}
\newcommand{\enthn}{\end{thm}}
\newcommand{\bepf}{\begin{proof}}
\newcommand{\enpf}{\end{proof}}
\newcommand{\belen}{\begin{lem}}
\newcommand{\enlen}{\end{lem}}
\newcommand{\becon}{\begin{cor}}
\newcommand{\encon}{\end{cor}}

\renewcommand{\vec}[1]{\bs{\mathrm{#1}}}

\newcommand{\basicspace}{{\cal S}}
\newcommand{\X}{\basicspace}

\newcommand{\supr}[1]{^{(#1)}}

\newcommand{\seq}[3]{#1_{#2}\ldots#1_{#3}}
\newcommand{\sseq}[3]{#1_{#2}^{#3}}  
\newcommand{\sumseq}[3]{\sum_{\sseq{#1}{#2}{#3}}}
\newcommand{\dsabs}[1]{\bigl| #1 \bigr|}
\newcommand{\labs}{\left| \vphantom{\sum_a^b} \right.}
\newcommand{\rabs}{\left. \vphantom{\sum_a^b} \right|}
\newcommand{\lparen}{\left( \vphantom{\sum} \right.}
\newcommand{\rparen}{\left. \vphantom{\sum} \right)}
\newcommand{\dsnrm}[1]{\bigl\Vert #1 \bigr\Vert}

\newcommand{\sd}{}
\newcommand{\scat}[2]{#1 \sd #2}

\newcommand{\rsub}[1]{_{#1}}  
\newcommand{\bsub}[1]{[#1]}   

\renewcommand{\phi}{\varphi}

\newcommand{\el}{\ell}

\newcommand{\calF}{\mathcal{F}}

\newcommand{\calL}{\mathcal{L}}

\newcommand{\calX}{\mathcal{X}}
\newcommand{\calY}{\mathcal{Y}}

\newcommand{\tha}{\theta}
\newcommand{\RR}{\mathbb{R}}
\newcommand{\NN}{\mathbb{N}}
\newcommand{\nn}[1]{[#1]}
\newcommand{\beq}{\begin{eqnarray*}}
\newcommand{\eeq}{\end{eqnarray*}}
\newcommand{\beqn}{\begin{eqnarray}}
\newcommand{\eeqn}{\end{eqnarray}}
\newcommand{\ben}{\begin{enumerate}}
\newcommand{\een}{\end{enumerate}}
\newcommand{\bit}{\begin{itemize}}
\newcommand{\eit}{\end{itemize}}
\newcommand{\paren}[1]{\left( #1 \right)}

\newcommand{\tlprn}[1]{\left\{ #1 \right\}}
\newcommand{\abs}[1]{\left| #1 \right|}
\newcommand{\nrm}[1]{\left\Vert #1 \right\Vert}
\newcommand{\floor}[1]{\ensuremath{\left\lfloor#1\right\rfloor}}
\newcommand{\gn}{\, | \,}
\newcommand{\ds}{\displaystyle}
\newcommand{\ts}{\textstyle}
\newcommand{\bs}{\boldsymbol}
\renewcommand{\th}{\ensuremath{^{\mathrm{th}}}~}

\def\blk{~\texttt{<LK>}~}
\def\elk{~\texttt{</LK>}~}
\newcommand{\lnote}[1]{\blk #1 \elk}
\newcommand{\hide}[1]{}
\newcommand{\oo}[1]{\frac{1}{#1}}
\newcommand{\defeq}{\doteq}

\title{Measure Concentration of Markov Tree Processes}
\author{Leonid Kontorovich\\
School of Computer Science\\ 
Carnegie Mellon University\\ 
Pittsburgh, PA 15213\\
USA \\
\url{lkontor@cs.cmu.edu}
}
\begin{document}
\maketitle
\abstract{We prove an apparently novel concentration of measure result
  for Markov tree processes. The bound we derive reduces to the known
  bounds for Markov processes when the tree is a chain, thus strictly
  generalizing the known Markov process concentration results. We employ
  several techniques of potential independent interest, especially for
  obtaining similar results for more general directed acyclic
  graphical models.
}

\newcommand{\lev}{\operatorname{lev}}
\newcommand{\rents}{\operatorname{parents}}
\newcommand{\kids}{\operatorname{children}}
\newcommand{\depth}{\operatorname{dep}}
\newcommand{\width}{\operatorname{wid}}
\newcommand{\dpath}{\vec \pi}

\newcommand{\tp}{\otimes}
\newcommand{\TP}{\bigotimes}
\newcommand{\s}{\sigma}
\newcommand{\bfs}{\bs{\sigma}}
\newcommand{\bfr}{\bs{\rho}}
\newcommand{\bft}{\bs{\tau}}

\section{Introduction}
\label{sec:intro}
An emerging paradigm for proving concentration results for nonproduct
measures is to quantify the dependence between the variables and state
the bounds in terms of that dependence (see \cite{kontram06} for an overview). 
A process (measure)
particularly 
amenable to
this approach is the Markov process.
Using different techniques, Marton (coupling method~\cite{marton96},
1996), Samson (log-Sobolev inequality~\cite{samson00},
  2000) and Kontorovich and Ramanan (martingale
  differences~\cite{kontram06}, 2006) have obtained qualitatively
  similar concentration of measure results for Markov processes.
One natural generalization of the Markov process is the {\it hidden
  Markov process}; we proved a concentration result for 
this class
in~\cite{kont06}. A different way to generalize the Markov process is
  via the {\it Markov tree process}, which we address in 
the present
paper.

If $(\X^n,d)$ is a metric space and 
$(X_i)_{1\leq i\leq n}$,
$X_i\in\X$ is a random process, a measure concentration result (for
the purposes of this paper) 
is an inequality
stating that for any 1-Lipschitz 
(with respect to $d$)
function $f:\X^n\to\RR$, we have
\beqn
\pr{\abs{f(X)-\mexp f(X)}>t} &\leq& 
2\exp(-Kt^2),
\eeqn
where $K$ may depend on $n$ but not on $f$.\footnote{See
\cite{ledoux01} for a much more general notion of concentration.}

\hide{
Recently several general
techniques have been developed for
proving concentration results for nonproduct measures~\cite{marton96,samson00}.
Let $(\X,\calF)$ be a
Borel-measurable space,
and consider the 
probability space
$(\X^n,\calF^n,\mu)$ 
with the associated random process $(X_i)_{1\leq i\leq n}$,
$X_i\in\X$. Suppose further that $\X^n$ is equipped with a metric $d$.
For our purposes, a concentration of measure result is an inequality
stating that for any Lipschitz 
(with respect to $d$)
function $f:\X^n\to\RR$, we have
\beqn
\pr{\abs{f(X)-\mexp f(X)}>t} &\leq& \alpha(t),
\eeqn
where $\alpha(t)$ is rapidly decaying to 0
as $t$ gets large.
}
The 
quantity 
$\bar\h_{ij}$, defined below,
has proved useful for obtaining concentration
results. For $1\leq i<j\leq n$, $y\in\X^{i-1}$ and $w\in\X$, let
$$\calL(\sseq{X}{j}{n}\gn \sseq{X}{1}{i-1}=y,X_i=w)
$$ be the law of $\sseq{X}{j}{n}$ conditioned on 
$\sseq{X}{1}{i-1}=y$ and $X_i=w$. Define
\beqn
\label{eq:hdef}
\h_{ij}(y,w,w') &=&
\TV{
\calL(\sseq{X}{j}{n}\gn\sseq{X}{1}{i-1}=y,X_i=w)-
\calL(\sseq{X}{j}{n}\gn\sseq{X}{1}{i-1}=y,X_i=w')}
\eeqn
and
\beq
\bar\h_{ij} &=&
\sup_{y\in\X^{i-1}} 
\sup_{w,w'\in\X} 
\h_{ij}(y,w,w')
\eeq
where $\TV{\cdot}$ is the total variation norm (see
\S\ref{sec:notconv} to clarify notation).

Let $\Gamma$ and $\Delta$ be upper-triangular $n\times n$ matrices,
with $\Gamma_{ii}=\Delta_{ii}=1$ and
\beq
\Gamma_{ij} = \sqrt{\bar\h_{ij}},
\qquad
\Delta_{ij} = \bar\h_{ij}
\eeq
for $1\leq i<j\leq n$.

For the case where $\X=[0,1]$
and
$d$ is the Euclidean metric on $\RR^n$, 
Samson~\cite{samson00} showed that
if
$f:[0,1]^n\to\RR$ is convex and Lipschitz with $\Lip{f}\leq1$, then
\beqn
\label{eq:samson}
\pr{\abs{f(X)-\mexp f(X)}>t} &\leq& 2\exp\paren{-\frac{t^2}{2\nrm{\Gamma}_2^2}}
\eeqn
where $\nrm{\Gamma}_2$ is the $\el_2$ operator norm of the matrix
$\Gamma$; Marton~\cite{marton98} has a comparable result.

For the case where $\X$ is countable and $d$ is the (normalized)
Hamming metric on $\X^n$,
$$ d(x,y) = \oo n\sum_{i=1}^n \pred{x_i\neq y_i},$$
Kontorovich and Ramanan~\cite{kontram06} showed that
if
$f:\X^n\to\RR$ is Lipschitz with $\Lip{f}\leq1$, then
\beqn
\label{eq:kontram}
\pr{\abs{f(X)-\mexp f(X)}>t} &\leq& 2\exp\paren{-\frac{nt^2}{2\nrm{\Delta}_\infty^2}}
\eeqn
where $\nrm{\Delta}_\infty$ is the $\el_\infty$ operator norm of the
matrix $\Delta$, also given by
\beqn
\nrm{\Delta}_\infty
&=&
\label{eq:infnorm}
\max_{1\leq i< n}(
1 + \bar\h_{i,i+1} + \ldots +\bar\h_{i,n}
).
\eeqn
This 
leads to
a strengthening of the Markov measure concentration result in
Marton~\cite{marton96}. 

The sharpest currently known Markov measure concentration results 
were
obtained in~\cite{kontram06}
and~\cite{samson00},
in terms of the contraction
coefficients 
$(\tha_i)_{1\leq i<n}$
of the Markov process:
\beqn
\label{eq:markbd}
\bar\h_{ij} &\leq& \tha_i\tha_{i+1}\cdots\tha_{j-1}.
\eeqn

In this 
paper,
we prove a bound on $\bar\h_{ij}$ in terms of the contraction
coefficients of the Markov tree process
(Theorem~\ref{thm:main}). This bound is cumbersome to
state without preliminary definitions, but it reduces to
(\ref{eq:markbd}) in the case where the Markov tree is a chain.

\section{Bounding 
$\bar\h_{ij}$ 
for Markov tree processes}
\subsection{Notational 
preliminaries
}
\label{sec:notconv}
Random variables are capitalized ($X$), specified state sequences
are written in lowercase ($x$), the shorthand
$\sseq{X}{i}{j}\equiv\seq{X}{i}{j}$ is used for all sequences, and
the 
concatenation
of the
sequences $x$ and $y$ is denoted by $\scat{x}{y}$, as in
$\scat{\sseq{x}{i}{j}}{\sseq{x}{j+1}{k}}=\sseq{x}{i}{k}$.
Another way to index collections of variables is by subset: if
$I=\{i_1, i_2,\ldots, i_m\}$ then we write 
$x_I\equiv x[I]\defeq \{x_{i_1},x_{i_2},\ldots,x_{i_m}\}$;
we will write $x_I$ and $x[I]$ interchangeably, as dictated by
convenience.
To avoid cumbersome subscripts, we will also occasionally use the
bracket notation for vector components. 
Thus, 
$\vec u\in \RR^{\X^I}$, then
$$ \vec u_{x_I}\equiv \vec u_{x[I]} \equiv \vec u[x[I]]
\defeq \vec u_{(x_{i_1},x_{i_2},\ldots,x_{i_m})} \in \RR$$
for each $x[I]\in\X^I$.
A similar bracket notation will apply for matrices.

We will use $\abs{\cdot}$ to denote set cardinalities.
Sums will range over the entire space of the summation variable;
thus
$\ds\sum_{\sseq{x}{i}{j}}f(\sseq{x}{i}{j})$ stands for
$\ds\sum_{\sseq{x}{i}{j}\in\X^{j-i+1}}f(\sseq{x}{i}{j})$,
and 
$\ds\sum_{x[I]}f(x[I])$ is shorthand for
$\ds\sum_{x[I]\in\X^I}f(x[I])$.

The probability operator $\pr{\cdot}$ is defined 
with respect 
the measure space specified in context.

We will write $\nn{n}$ for the set $\{1,\ldots,n\}$. Anytime
$\nrm{\cdot}$ appears without a subscript, it will always denote the
{\bf total variation} norm $\TV{\cdot}$,
which we define here, for any 
signed measure 
$\tau$
on a countable set
$\calX$, by
\beqn
\label{eq:tv}
 \TV{\tau} 
\defeq
{\ts\oo2}\sum_{x\in\calX}\abs{\tau(x)}.
\eeqn
\hide{
 Regarding the latter, we
recall that
if
$\tau$ is a signed, balanced measure on a countable set $\calX$ (i.e.,
$\tau(\calX)=\sum_{x\in\calX}\tau(x)=0$), then
\beqn
\label{eq:tv}
 \TV{\tau} = {\ts\oo2}\nrm{\tau}_1 \equiv
{\ts\oo2}\sum_{x\in\calX}\abs{\tau(x)}.
\eeqn
}

If $G=(V,E)$ is a graph,
we will 
frequently
abuse notation and write $u\in G$ instead of $u\in
V$, blurring the distinction between a graph and its vertex set. This
notation will carry over to set-theoretic operations ($G=G_1\cap G_2$)
and indexing of variables (e.g., $X_G$).

Unless we will need to refer explicitly to a $\sigma$-algebra, we will
suppress it in the probability space notation, using less 
rigorous
formulations, such as ``Let $\mu$ be a measure on
$\X^n$''. Furthermore, to avoid the technical but inessential
complications associated with infinite sets, we will take $\X$ to be
finite in this paper, noting only that the bounds carry over unchanged to
the countable case (as done in~\cite{kontram06} and~\cite{kont06}). To
extend the results to the continuous case, some 
mild
measure-theoretic assumptions are needed (see~\cite{marton98}).

\subsection{Definition of Markov tree process}
\label{sec:mtpdef}
\subsubsection{Graph-theoretic preliminaries}
\label{sec:graphtheo}
Consider a directed acyclic graph 
$G=(V,E)$, and define a
partial order 
$\prec_G$
on $G$ by the transitive closure of the relation
$$u\prec_G v 
\qquad\text{if}\qquad
(u,v)\in E.$$
We define the {\bf parents} and {\bf children} of $v\in V$ in the
natural way:
$$ \rents(v) = \{u\in V : (u,v)\in E\} $$
and
$$ \kids(v) = \{w\in V : (v,w)\in E\}. $$

If $G$ is connected and each $v\in V$ has at most one parent, $G$ is called a
{\bf(directed) tree}. In a
tree, whenever $u\prec_G v$ there is a unique directed path from $u$ to
$v$.
A tree $T$ always has a unique minimal (w.r.t.
$\prec_T$)
element $r_0\in V$, called its {\bf root}. Thus, for every $v\in V$
there is a unique directed path 
$r_0 \prec_T r_1 \prec_T \ldots \prec_T r_d = v$; define
the {\bf depth} of $v$,
$\depth_T(v)=d$, to be the length (i.e., number of edges) of this path.
Note that $\depth_T(r_0)=0$. We define the depth of the tree by
$\depth(T)=\sup_{v\in T}\depth_T(v)$.

For $d=0,1,\ldots$ define
the $d$\th {\bf level} of the tree $T$ by
$$ \lev_T(d) = \{v \in V : \depth_T(v)=d\};$$
note that the levels induce a disjoint partition on $V$:
$$ V = \bigcup_{d=1}^{\depth(T)} \lev_T(d).$$
We define the {\bf width} of a tree as the greatest number of nodes in
any level:
\beqn
\label{eq:wid}
 \width(T) = \sup_{1\leq d\leq\depth(T)} \abs{\lev_T(d)}.
\eeqn

We will consistently take $|V|=n$ for finite $V$.
An ordering 
$J:V\to\NN$
of the nodes is 
said to be
{\bf breadth-first} if 
\beqn
 \depth_T(u) < \depth_T(v) \Longrightarrow J(u) < J(v) .
\eeqn
Since every directed tree
$T=(V,E)$ has some breadth-first ordering,\footnote{
One can easily construct a breadth-first ordering on a given tree by
ordering the nodes arbitrarily within each level and listing the
levels in ascending order: $\lev_T(1),\lev_T(2),\ldots$.
}
we shall henceforth blur
the distinction between $v\in V$ and $J(v)$, 
simply taking 
$V=\nn{n}$
(or $V=\NN$)
and assuming that $\depth_T(u) < \depth_T(v) \Rightarrow u < v$ holds. This
will allow us to write $\X^V$ simply as $\X^n$
for any set $\X$.

Note that we have two orders on $V$: the partial order $\prec_T$,
induced by the tree topology, and the total order $<$, given by the
breadth-first enumeration. Observe that $i\prec_T j$ implies $i<j$ but
not 
vice versa.

If $T=(V,E)$ is a tree and $u\in V$, we define the {\bf subtree} 
induced by $u$,
$T_u=(V_u,E_u)$
by $V_u = \{v\in V: u\preceq_T v\}$, 
$E_u = \{(v,w) \in E: v,w\in V_u\}$.

\subsubsection{Markov tree measure}
\label{sec:MTmeas}
If $\X$ is a finite set,
a
{\bf Markov tree measure} $\mu$ is defined on
$\X^n$ 
by a tree $T=(V,E)$
and transition kernels $p_0$, 
$\tlprn{p_{ij}(\cdot\gn\cdot)}_{(i,j)\in E}$.
Continuing our convention in \S\ref{sec:graphtheo}, we have a
breadth-first order $<$ and the total order $\prec_T$ on $V$, and take
$V=\{1,\ldots,n\}$. Together, the topology of $T$ and the transition
kernels determine the measure $\mu$ on $\X^n$:
\beqn
\label{eq:mtmeas}
\mu(x) = p_0(x_1)\prod_{(i,j)\in E} p_{ij}(x_j\gn x_i).
\eeqn
A measure
on 
$\X^n$ 
satisfying
(\ref{eq:mtmeas}) for some $T$ and $\{p_{ij}\}$
is said to be {\bf compatible} with tree $T$;
a measure is
a Markov tree measure if it is compatible with some tree.

Suppose $\X$ is a finite set and $(X_i)_{i\in\NN}$, $X_i\in\X$ is a
random process 
defined on 
$(\X^\NN,\P)$. 
If for each $n>0$
there is a tree $T\supr n=(\nn{n},E\supr n)$ 
and a Markov tree measure $\mu_n$
compatible with $T\supr n$
such that
for all $x\in\X^n$ we have
$$ \pr{\sseq{X}{1}{n} = x} = \mu_n(x)$$
then we call $X$ a {\bf Markov tree process}. 
The trees $\{T\supr n\}$ are easily seen to be consistent in the sense
that $T\supr n$ is an induced subgraph of $T\supr{n+1}$.
So corresponding 
to 
any Markov tree process is the unique infinite tree
$T=(\NN,E)$. 
The uniqueness of $T$ 
is easy to see,
since for
$v>1$, the parent of $v$ is the 
smallest
$u\in \NN$ such that
$$ \pr{X_v = x_v \gn \sseq{X}{1}{u}\;=\;\sseq{x}{1}{u}} = \pr{X_v = x_v \gn X_u=x_u};$$
thus $\P$ determines the topology of $T$. 

It is straightforward to
verify 
that 
a Markov tree process $\{X_v\}_{v\in T}$
compatible with tree $T$
has the following {\bf Markov property}: if $v$ and $v'$ are children
of $u$ in $T$, then
$$ \pr{X\rsub{T_v}=x,X\rsub{T_{v'}}=x'\gn X_{u}=y}
= \pr{X\rsub{T_v}=x\gn X_{u}=y}
\pr{X\rsub{T_{v'}}=x'\gn X_{u}=y}.$$
In other words, the subtrees induced by the children are conditionally
independent given the parent; this follows directly from the
definition of the Markov tree measure in (\ref{eq:mtmeas}).

\subsection{Statement of result}
\bethn
\label{thm:main}
Let $\X$ be a finite set and
let $(X_i)_{1\leq i\leq n}$, $X_i\in\X$  be a Markov tree process, defined by a
tree $T=(V,E)$ and transition kernels $p_0$, 
$\tlprn{p_{uv}(\cdot\gn\cdot)}_{(u,v)\in E}$.
Define the $(u,v)$- {\bf 
contraction coefficient} $\tha_{uv}$ by
\beqn
\label{eq:thadef}
\tha_{uv} &=&
\max_{y,y'\in\X}
\TV{
p_{uv}(\cdot\gn y)-p_{uv}(\cdot\gn y')
}.
\eeqn
Suppose 
$\max_{(u,v)\in E} \tha_{uv}\leq \tha <1 $
for some $\tha$ 
and
$\width(T)\leq L$
\hide{
$$\max
\tlprn{
\abs{\lev(d)} 
:
d=\depth(u), u\in\nn{n}
}
\;\leq\; L
$$
for some $L$
(that is, any level of $T$ contains at most $L$ nodes). 
}. 
Then for the
Markov tree process $X$ we have 
\beqn
\label{eq:mainbd}
\bar\h_{ij} &\leq& \paren{1-(1-\tha)^L}^{\floor{(j-i)/L}}
\eeqn
for $1\leq i<j\leq n$.
\enthn

\hide{
\begin{rem}
\lnote{remove/modify}
Modulo measurability issues,
a hidden Markov process may be defined on continuous hidden and
observed state spaces; the definition of
$\bar\h_{ij}$ is unchanged. For convenience, the proof of 
Theorem~\ref{thm:main} is given for the countable case, but can
straightforwardly
be extended to the continuous one.
\end{rem}

\lnote{actual bound is more refined; compare to other results
; effective $\tha$
}
}

To cast 
(\ref{eq:mainbd}) in more usable form, we first note that for $L\in\NN$
and $k\in\NN$, if $k\geq L$ then
\beqn
\label{eq:flrbd}
 \floor{\frac{k}{L}} \geq \frac{k}{2L-1}
\eeqn
(we omit the elementary number-theoretic proof). Using
(\ref{eq:flrbd}), we have
\beqn
\label{eq:thatil}
\bar\h_{ij} &\leq& \tilde\tha^{j-i},
\qquad\text{for } j\geq i+L
\eeqn
where
$$ \tilde\tha = (1-(1-\tha)^L)^{
1/(
{2L-1}
)
}.$$

The bounds in (\ref{eq:samson}) and (\ref{eq:kontram}) are for
different metric spaces and therefore not readily comparable
(the result in (\ref{eq:samson}) has the additional convexity
assumption; see~\cite{kont06-metric-mix} for a discussion).\hide{
In the
special
case where the underlying Markov process is contracting, i.e.,
$\tha_i\leq\tha<1$ for $1\leq i<n$, 
Theorem~\ref{thm:main} yields
\beq
\bar\h_{ij} &\leq& \tha^{j-i}.
\eeq
In this case, 
}
For the case where (\ref{eq:thatil}) holds,
Samson's bound \cite{samson00} yields
\beqn
\label{eq:gamapprox}
\nrm{\Gamma}_2 &\lesssim& \oo{1-{\tilde\tha}^{\oo2}},
\eeqn
and the 
approximation
\beqn
\label{eq:delapprox}
\nrm{\Delta}_\infty &\lesssim& \sum_{k=0}^\infty\tilde\tha^k=\oo{1-\tilde\tha}
\eeqn
holds trivially via (\ref{eq:infnorm}).\hide{
\footnote{
The statement is approximate
because 
(\ref{eq:thatil}) does not hold for all $j>i$ but only starting with
$j\geq i+L$. The 
difference between
$\paren{1-(1-\tha)^L}^{\floor{(j-i)/L}}=1$
and
$\tilde\tha^{j-i}$
for $i<j<i+L$ 
is
at most
$1-\tilde\tha^{L-1}$
and affects only a fixed finite number
($L-1$)
 of entries in each row of
$\Gamma$
and $\Delta$. Since $\nrm{\cdot}_2$ and 
$\nrm{\cdot}_\infty$ are continuous functionals, we are justified in
claiming the approximate bound, which may be quantified if an
application calls for it.
The statements in (\ref{eq:gamapprox}) and (\ref{eq:delapprox}) are
only meant to convey an order of magnitude.
}}
In the (degenerate) case where
the Markov tree is a chain, we have $L=1$ and therefore
$\tilde\tha=\tha$; thus we recover the Markov chain concentration
results in~\cite{kontram06,marton96,samson00}
and the approximations in (\ref{eq:gamapprox},\ref{eq:delapprox})
become precise inequalities.
\begin{rem}
The bounds in 
(\ref{eq:gamapprox}) and (\ref{eq:delapprox}) are
approximate
because 
(\ref{eq:thatil}) does not hold for all $j>i$ but only starting with
$j\geq i+L$. The 
difference between
$\paren{1-(1-\tha)^L}^{\floor{(j-i)/L}}=1$
and
$\tilde\tha^{j-i}$
for $i<j<i+L$ 
is
at most
$1-\tilde\tha^{L-1}$
and affects only a fixed finite number
($L-1$)
 of entries in each row of
$\Gamma$
and $\Delta$. Since $\nrm{\cdot}_2$ and 
$\nrm{\cdot}_\infty$ are continuous functionals, we are justified in
claiming the approximate bound, which may be quantified if an
application calls for it.
The statements in (\ref{eq:gamapprox}) and (\ref{eq:delapprox}) are
only meant to convey an order of magnitude.
\end{rem}

\subsection{Proof of 
main result
}
\label{sec:proofmain}
\hide{
Let us first quote a simple lemma from~\cite{kontram06}, which 
has probably been known for quite some time:
}
The proof of Theorem~\ref{thm:main} is combination of elementary graph
theory and tensor algebra. We start with a graph-theoretic lemma:
\belen
\label{lem:j0}
Let $T=(\nn{n},E)$ be a tree and fix $1\leq i<j\leq n$. Suppose 
$(X_i)_{1\leq i\leq n}$ is a Markov tree process 
whose
law
$\P$
on $\X^n$ 
is
compatible with $T$ (in the sense of
\S\ref{sec:MTmeas}). Define the set
$$ 
\hide{
T_{ij} = \{
 k\in T_i :
k\geq j\}
}
T_i^j = T_i \cap \{j,j+1,\ldots, n\}
, $$
consisting of those nodes in the subtree $T_i$ whose
breadth-first numbering does not precede $j$. Then, for $y\in\X^{i-1}$
and $w,w'\in\X$, we have
\beqn
\h_{ij}(y,w,w') &=&
\left\{
\begin{array}{ll}
0 & \quad T_i^j = \emptyset \\
\h_{ij_0}(y,w,w') & \quad \mbox{otherwise},
\end{array}
\right.
\eeqn
where $j_0$ is the minimum (with respect to $<$) element of $T_i^j$.
\enlen
\begin{rem}
This lemma tells us that when computing $\h_{ij}$ it is sufficient to
restrict our attention to the subtree induced by $i$.
\end{rem}
\bepf
The case $j\in T_i$ implies $j_0=j$ and is trivial; thus we assume
$j\notin T_i$. In this case, the subtrees $T_i$ and $T_j$ are
disjoint. Putting $\bar T_i=T_i\setminus\{i\}$, 
we have by the Markov 
property, 
\beq
\pr{X\rsub{\bar T_i}=x\rsub{\bar T_i},X\rsub{T_j}=x\rsub{T_j}\gn 
\sseq{X}{1}{i} = \scat{y}{w}}
&=& 
\pr{X\rsub{\bar T_i}=x\rsub{\bar T_i}\gn X_i=w}
  \pr{X\rsub{T_j}=x\rsub{T_j}\gn \sseq{X}{1}{i-1} = y}.
\eeq

Then
from (\ref{eq:hdef}) and (\ref{eq:tv}), and by
marginalizing out the $X_{T_j}$,
we have
\beq
\h_{ij}(y,w,w') &=&
{\ts\oo2}
\sumseq{x}{j}{n}
\abs{
\pr{\sseq{X}{j}{n}=\sseq{x}{j}{n}\gn\sseq{X}{1}{i}=\scat{y}{w}}
-    
\pr{\sseq{X}{j}{n}=\sseq{x}{j}{n}\gn\sseq{X}{1}{i}=\scat{y}{w'}}
}\\
&=&
{\ts\oo2}
\sum_{x\rsub{T_i^j}}
\abs{
\pr{X\rsub{T_i^j}=x\rsub{T_i^j}\gn X_i=w}
-
\pr{X\rsub{T_i^j}=x\rsub{T_i^j}\gn X_i=w'}
}.
\eeq
If $T_i^j=\emptyset$ then obviously $\h_{ij}=0$; otherwise,
$\h_{ij}=\h_{ij_0}$, since $j_0$ is the 
``first''
element of $T_i^j$.
\hide{
This leaves two possible cases:
\bit
\item[(\one)] there is a (smallest) $j_0\in T_i$, $j_0\geq j$
\item[(\two)] for each $k\in T_i$, $k<j$.
\eit
}
\enpf

Next we develop some basic results for tensor norms; recall that
unless specified otherwise, the
norm used in this paper is the total variation norm defined in (\ref{eq:tv}).
If 
$\vec A$
is an $M\times N$
column-stochastic 
matrix:
($\vec A_{ij}\geq0$ for $1\leq i\leq M$, $1\leq j\leq N$ and
$\sum_{i=1}^M \vec A_{ij}=1$ for all $1\leq j\leq N$) and $\vec u\in \RR^N$
is {\it balanced} in the sense that $\sum_{j=1}^N \vec u_j=0$, we
have, by 
the 
Markov contraction lemma
(\cite{kontram06}, Lemma B.1),
\beqn
\label{eq:contr}
\nrm{\vec{Au}} &\leq& \nrm{\vec A}\nrm{\vec u},
\eeqn
where
\beqn
\label{eq:matnorm}
\nrm{\vec A} &=& \max_{1\leq j,j'\leq N} \nrm{\vec A_{*,j}-\vec A_{*,j'}},
\eeqn
and 
$\vec A\rsub{*,j}
\equiv
\vec A\bsub{\cdot,j}
$ 
denotes the $j$\th column of $\vec A$. An immediate
consequence
of (\ref{eq:contr}) is that $\nrm{\cdot}$ 
satisfies
\beqn
\label{eq:AB}
\nrm{\vec{AB}} &\leq& \nrm{\vec A}\nrm{\vec B}
\eeqn
for column-stochastic matrices
$\vec A\in\RR^{M\times N}$ and $\vec B\in\RR^{N\times P}$.
\begin{rem}
\label{rem:stochnorm}
Note that 
if $\vec A$ is a column-stochastic matrix
then
$\nrm{\vec A}\leq1$, 
and if additionally $\vec u$ is balanced then 
$\vec{Au}$ is also balanced.
\end{rem}

If $\vec u\in\RR^M$ and $\vec v\in\RR^N$, define their tensor product
$\vec w = \vec v \tp \vec u$
by
\beq 
\vec w_{(i,j)} &=& \vec u_i \vec v_j,
\eeq
where the notation $(\vec v \tp \vec u)_{(i,j)}$ is used to
distinguish the 2-tensor $\vec w$ from an $M\times N$ matrix. The
tensor $\vec w$ is a vector in $\RR^{MN}$ indexed by
pairs $(i,j)\in\nn{M}\times\nn{N}$; its norm is naturally defined to
be
\beqn
\label{eq:tensnorm}
\nrm{\vec w} = {\ts\oo2}\sum_{(i,j)\in\nn{M}\times\nn{N}}\abs{\vec w_{(i,j)}}.
\eeqn

The following 
``tensorizing'' lemma
will play a key role in deriving our bound
(we suppress the boldfaced vector notation for readability):
\renewcommand{\X}{\calX}
\newcommand{\Y}{\calY}
\belen
\label{lem:tvineq}
Consider two
finite sets $\X,\Y$, with probability measures $p,p'$ on $\X$ and
$q,q'$ on $\Y$.
Then 
\beqn
\nrm{p\tp q - p'\tp q'} &\leq &
\nrm{p-p'} + \nrm{q-q'}
- \nrm{p-p'}\nrm{q-q'}.
\eeqn
\enlen
\begin{rem}
Note that $p\tp q$ is a 2-tensor in $\RR^{\X\times\Y}$ and a
probability measure on $\X\times\Y$.
\end{rem}
\bepf
\hide{
The claim amounts to showing that for $p,p',q,q'$ as above, we have 
\beqn
\nonumber
\sum_{x\in\X,y\in\Y}\abs{p_x q_y - p'_x q'_y}
&\leq&
\sum_{x\in\X}\abs{p_x - p'_x} +
\nrm{q-q'}
-
\paren{\sum_{x}\abs{p_x - p'_x}}
\nrm{q-q'}
.\\
\label{eq:tvtensclaim}
\eeqn
}
Fix $q,q'$ and define the function
$$ F(u,v) = 
\sum_{x\in\X}\abs{u_x - v_x} +
\nrm{q-q'}
\paren{2
-
\sum_{x\in\X}\abs{u_x - v_x}
}
-
\sum_{x\in\X,y\in\Y}\abs{u_x q_y - v_x q'_y}
$$
over the convex polytope
$U\subset \RR^{\X}\times\RR^{\X}$,
$$ U = \tlprn{(u,v): 
u_x,v_x\geq0,\sum u_x=\sum v_x=1};$$
note that proving 
the claim
is equivalent to showing that
$F\geq0$ on $U$.

For any $\bfs\in\{-1,+1\}^{\X}$, let
\beq
U_{\bfs} = \{(u,v)\in U : \sgn( u_x- v_x)=\s_x\};
\eeq
note that $U_{\bfs}$ is a convex polytope
and that $U = \bigcup_{\bfs\in\{-1,+1\}^\X} U_{\bfs}$.\footnote{
We define $\sgn(z) = \pred{z\geq0} - \pred{z<0}$. Note that 
the constraint 
$\sum_{x\in\X}  u_x =
 \sum_{x\in\X}  v_x = 1$
forces
$ U_{\bfs} = \{(u,v)\in U :  u_x= v_x\}  $
when
$\s_x=+1$ for all $x\in\X$
and $U_{\bfs} = \emptyset$
when
$\s_x=-1$ for all $x\in\X$.
Both of these cases are trivial.
}

Pick an arbitrary
$\bft\in\{-1,+1\}^{\X\times\Y}$ 
and
define
\beq
F_{\bfs}(u,v) = 
\sum_{x}\s_x(u_x - v_x) +
\nrm{q-q'}
\paren{2-\sum_{x}\s_x(u_x - v_x)}
-
\sum_{x,y}\tau_{xy}(u_x q_y - v_x q'_y)
\eeq
over $U_{\bfs}$.
Since 
$\s_x(u_x-v_x)=\abs{u_x-v_x}$
and 
%
$\tau_{xy}$ can be chosen 
(for any given $u,v,q,q'$)
so that
$\tau_{xy}(u_x q_y - v_x q'_y)=\abs{u_x q_y - v_x q'_y}$,
the claim that 
$F \geq0$ on $U$
will follow if we can show that
$F_{\bfs} \geq0$ on $U_{\bfs}$.

Observe that $F_{\bfs}$ is affine in its arguments $(u,v)$ and recall that an
affine function achieves its extreme values on the extreme points of a
convex domain. 
Thus to verify that $F_{\bfs} \geq0$ on $U_{\bfs}$, we need only check the
value of $F_{\bfs}$ on the extreme points of $U_{\bfs}$.
The extreme points of $U_{\bfs}$ are pairs $(u,v)$ such that,
for some $x',x''\in\X$,
$u=\delta(x')$
and
$v=\delta(x'')$,
where $\delta(x_0)\in\RR^\X$ is given by $[\delta(x_0)]_x = \pred{x=x_0}$.
Let $(\hat u,\hat v)$ be an extreme point of $U_{\bfs}$. The case 
$\hat u=\hat v$ 
is trivial,
so assume 
$\hat u\neq \hat v$.
 In this case,
$
\sum_{x\in\X} \s_x(\hat u_x-\hat v_x) = 2
$
and
\beq
\hide{
\sum_{x\in\X,y\in\Y}\tau_{xy}(\hat u_x q_y - \hat v_x q'_y)
&=&
\sum_{x\in\X,y\in\Y}\tau_{xy}([\delta(x')]_x q_y - [\delta(x'')]_x q'_y)\\
&=&
\sum_{y} \tau_{x'y}q_y - \sum_y \tau_{x''y}q'_y,
}
\dsabs{\sum_{x\in\X,y\in\Y}\tau_{xy}(\hat u_x q_y - \hat v_x q'_y)}
&\leq&
\sum_{x\in\X,y\in\Y}\abs{\hat u_x q_y - \hat v_x q'_y}\\
&\leq& 2.
\eeq

This shows that 
$F_{\bfs}\geq 0$
on $U_{\bfs}$ and completes the proof.
\enpf
\renewcommand{\X}{\basicspace}

To develop a convenient tensor notation, we will fix the index set
$V=\{1,\ldots,n\}$. 
For $I\subset V$, a tensor indexed by $I$ is a 
vector $\vec u\in\RR^{\X^I}$.
A special case of such an 
$I$-tensor is the product
$\vec u = \TP_{i\in I} \vec v\supr i$, where $\vec v\supr i\in\RR^\X$ and
%
\beq
\vec u\bsub{x\rsub{I}} &=&
\prod_{i\in I} 
\vec v\supr{i}\bsub{x_i}
\eeq
for each $x\rsub{I}\in{\X^I}$.
To gain more familiarity with the notation, let us write the total
variation norm of an $I$-tensor:
\beqn
\label{eq:Itensnorm}
\nrm{\vec u} &=&
{\ts\oo2}\sum_{x\rsub{I}\in\X^I} \abs{\vec u\bsub{x\rsub{I}}}.
\eeqn
In order to extend
Lemma~\ref{lem:tvineq} 
to product tensors, we will need to define 
the
function
$\inex_k:\RR^k\to\RR$ and state some of its properties:
\belen
\label{lem:inex}
Define $\inex_k:\RR^k\to\RR$ recursively as $\inex_1(x)=x$ and
\beqn
\label{eq:inex}
\inex_{k+1}(x_1,x_2,\ldots,x_{k+1})=
x_{k+1}+(1-x_{k+1})\inex_{k}(x_1,x_2,\ldots,x_{k}).
\eeqn
Then
\bit
\item[\iia] $\inex_k$ is symmetric in its $k$
arguments, so it is well-defined as a mapping 
$$\inex:\{x_i:1\leq i\leq k\}\mapsto\RR$$ from finite real sets to the reals
\item[\iib] $\inex_k$ takes $[0,1]^k$ to $[0,1]$ and is monotonically
      increasing in each argument on $[0,1]^k$
\item[\iic] If $B\subset C\subset[0,1]$ 
are finite sets
then $\inex(B)\leq\inex(C)$
\item[\iid] $\inex_k(x,x,\ldots,x)=1-(1-x)^k$
\item[\iie] if $B$ is finite and $1\in B\subset[0,1]$ then $\inex(B)=1$.
\item[\iif] if $B\subset[0,1]$
is a finite set then
 $\inex(B)\leq\sum_{x\in B} x$.
\eit
\enlen
\begin{rem}
In light of~\iia, we will use the notation
$\inex_k(x_1,x_2,\ldots,x_{k})$ and $\inex(\{x_i:1\leq i\leq k\})$
interchangeably, 
as dictated by convenience.
\end{rem}
\bepf
Claims~\iia,~\iib,~\iie,~\iif~are straightforward to verify from the
recursive 
definition of $\inex$ and induction. Claim~\iic~follows from~\iib~since
$$
\inex_{k+1}(x_1,x_2,\ldots,x_{k},0)=\inex_{k}(x_1,x_2,\ldots,x_{k})$$
and~\iid~is easily derived from the binomial expansion of $(1-x)^k$.
\enpf

The function $\inex_k$ is the natural generalization of 
$\inex_2(x_1,x_2)=x_1+x_2-x_1x_2$
to $k$ variables, and it is what we need 
for the analogue of Lemma~\ref{lem:tvineq} 
for a  product of $k$ tensors:
\becon
\label{cor:tp}
Let
$\{\vec u\supr i\}
_{i\in I}
$ 
and 
$\{\vec v\supr i\}
_{i\in I}
$
be
two sets of tensors
and
assume 
that each of 
$\vec u\supr i,\vec v\supr i$ is a probability measure on $\X$. Then
we have
\beqn
\nrm{\TP_{i\in I}\vec u\supr i - \TP_{i\in I}\vec v\supr i}
&\leq&
\inex
\tlprn{
\dsnrm{\vec u\supr i-\vec v\supr i}: i\in I
}
\hide{
\dsnrm{\vec u\supr1-\vec u\supr1},
\dsnrm{\vec u\supr2-\vec u\supr2},
\ldots,
\dsnrm{\vec u\supr{k}-\vec u\supr{k}})
}.
\eeqn
\hide{
where $\inex_k:\RR^k\to\RR$ is defined recursively as $\inex_1(x)=x$ and
\beqn
\label{eq:inex}
\inex_{k+1}(x_1,x_2,\ldots,x_{k+1})=
x_{k+1}+(1-x_{k+1})\inex_{k}(x_1,x_2,\ldots,x_{k}).
\eeqn
\encon
\begin{rem}
\label{rem:inex}
It is clear that the function $\inex_k$ is symmetric in its $k$
arguments, so it is well-defined as a mapping 
$$\inex:\{x_i:1\leq i\leq k\}\mapsto\RR$$ from finite real sets to the reals.
\end{rem}
}
\encon
\bepf
Pick an $i_0\in I$ and let 
$\vec p=\vec u\supr{i_0}$,
$\vec q=\vec v\supr{i_0}$,
$$\vec p' = \TP_{i_0\neq i\in I} \vec u\supr i,
\qquad
\vec q' = \TP_{i_0\neq i\in I} \vec v\supr i.$$
Apply Lemma~\ref{lem:tvineq} to 
$\nrm{\vec p\tp\vec q-\vec p'\tp\vec q'}$ and proceed by induction.
\enpf

Our final generalization concerns linear operators over
$I$-tensors. An $I,J$-matrix $\vec A$ 
has dimensions $|\X^J|\times|\X^I|$ and
takes an $I$-tensor $\vec u$ to
a $J$-tensor $\vec v$:
for each $y\rsub{J}\in\X^J$,
we have
\beqn
\label{eq:Au}
\vec{v}\bsub{y\rsub{J}} 
&=&
\sum_{x\rsub{I}\in\X^I} 
\vec A\bsub{y\rsub{J},x\rsub{I}} \vec u\bsub{x\rsub{I}},
\eeqn
which we write as $\vec{Au}=\vec{v}$.
If $\vec A$ is an $I,J$-matrix and $\vec B$ is a $J,K$-matrix, the
matrix product $\vec{BA}$ is defined analogously to (\ref{eq:Au}).

As a special case, an $I,J$-matrix might factorize as a tensor product
of $|\X|\times|\X|$ matrices 
$\vec A\supr{i,j}\in\RR^{\X\times\X}$.
We will write such a factorization in terms of a bipartite
graph\footnote{
Our notation for bipartite graphs is standard; it is equivalent to
$G=(I\cup J,E)$ where $I$ and $J$ are always assumed to be disjoint.
}
$G=(I+J,E)$, where $E\subset I\times J$ and the factors 
$\vec A\supr{i,j}$ are indexed by $(i,j)\in E$:
\beqn
\label{eq:matensprod}
 \vec A = \TP_{(i,j)\in E} \vec A\supr{i,j}, 
\eeqn
where
$$
\vec A\bsub{y\rsub{J},x\rsub{I}}
= 
\prod_{(i,j)\in E} \vec A\supr{i,j}_{y_j,x_i}
$$
for all $x\rsub{I}\in\X^I$ and $y\rsub{J}\in\X^J$.
The norm of an $I,J$-matrix is a natural generalization of the matrix
norm defined in (\ref{eq:matnorm}):
\beqn
\label{eq:ijmatnorm}
\nrm{\vec A} = \max_{x\rsub{I},x\rsub{I}'\in\X^I} 
\nrm{
\vec A\bsub{\cdot,x\rsub{I}} - 
\vec A\bsub{\cdot,x\rsub{I}'}
}
\eeqn
where 
$\vec u=\vec A\bsub{\cdot,x\rsub{I}}$
is the $J$-tensor given by
$$ \vec u\bsub{y\rsub{J}} = 
\vec A\bsub{y\rsub{J},x\rsub{I}};
$$
(\ref{eq:ijmatnorm}) is well-defined via the tensor norm in 
(\ref{eq:Itensnorm}).
Since $I,J$ matrices act on $I$-tensors by
ordinary matrix multiplication, 
$\nrm{\vec{Au}} \leq \nrm{\vec A}\nrm{\vec u}$
continues to hold
when $\vec A$ is a
column-stochastic
$I,J$-matrix and $\vec u$ is a
balanced
$I$-tensor; if, additionally, $\vec B$ is a column-stochastic $J,K$-matrix,
$\nrm{\vec{BA}} \leq \nrm{\vec B}\nrm{\vec A}$
also holds.
Likewise, since another way of writing
(\ref{eq:matensprod}) is
\beq
\vec A\bsub{\cdot,x\rsub{I}} = \TP_{(i,j)\in E}
\vec A\supr{i,j}
\bsub{\cdot,x_i}
,
\eeq
Corollary~\ref{cor:tp} extends to tensor products of matrices:
\belen
\label{lem:TP}
Fix index sets 
$I,J$ 
and a 
bipartite graph $(I+J,E)$. 
Let 
$\tlprn{\vec A\supr{i,j}}_{(i,j)\in E}$ be a collection of column-stochastic
$|\X|\times|\X|$ matrices, whose tensor product is the $I,J$ matrix
$$
 \vec A = 
\TP_{(i,j)\in E} \vec A\supr{i,j}.
$$
Then
\beq
\nrm{\vec A} &\leq&
\inex
\tlprn{
\dsnrm{\vec A\supr{i,j}}:(i,j)\in E
}
\hide{
(
\dsnrm{\vec A\supr{i_1,j_1}},
\dsnrm{\vec A\supr{i_2,j_2}},
\ldots,
\dsnrm{\vec A\supr{i_k,j_k}})
}.
\eeq
\enlen

We are now in a position to state the main technical lemma, from which
Theorem~\ref{thm:main} will follow straightforwardly:
\belen
\label{lem:maintech}
Let $\X$ be a finite set and
let $(X_i)_{1\leq i\leq n}$, $X_i\in\X$  be a Markov tree process, defined by a
tree $T=(V,E)$ and transition kernels $p_0$, 
$\tlprn{p_{uv}(\cdot\gn\cdot)}_{(u,v)\in E}$.
Let the $(u,v)$-contraction coefficient
$\tha_{uv}$
be as defined in (\ref{eq:thadef}).

Fix $1\leq i<j\leq n$ and let 
$j_0=j_0(i,j)$ be as defined in Lemma~\ref{lem:j0} (we are
assuming its existence, for otherwise $\bar\h_{ij}=0$).
Then we have
\beqn
\bar\h_{ij} &\leq&
\prod_{d=
\depth_T(i)+1
}^{
\depth_{T}(j_0)
}
\inex
\tlprn{
\tha_{uv} : v\in\lev_{T}(d)
}
\eeqn
where $\depth_{T}(\cdot)$ is defined in \S\ref{sec:graphtheo}.
\enlen
\newcommand{\sz}{\sseq{z}{i+1}{j-1}}
\newcommand{\sx}{\sseq{x}{j}{n}}
\bepf
For $y\in\X^{i-1}$ and $w,w'\in\X$, we have
\beqn
\h_{ij}(y,w,w') &=&
{\ts\oo2}
\sumseq{x}{j}{n}
\abs{
\pr{\sseq{X}{j}{n}=\sseq{x}{j}{n}\gn\sseq{X}{1}{i}=\scat{y}{w}}
-    
\pr{\sseq{X}{j}{n}=\sseq{x}{j}{n}\gn\sseq{X}{1}{i}=\scat{y}{w'}}
}\\
&=&
{\ts\oo2}
\sumseq{x}{j}{n}
\labs
\sumseq{z}{i+1}{j-1}
\lparen{
\pr{\sseq{X}{i+1}{n}=\sz \sx\gn\sseq{X}{1}{i}=\scat{y}{w}}
}\nonumber\\
&&
\qquad\qquad\quad
-    
\pr{\sseq{X}{i+1}{n}=[\sz\,\sx]\gn\sseq{X}{1}{i}=\scat{y}{w'}}
\rparen\rabs.
\eeqn
Let $T_i$ be the subtree induced by $i$ and
\beqn
\label{eq:ZC}
Z=T_i\cap\{i+1,\ldots,j_0-1\}
\qquad\text{and}\qquad
C = \{v\in T_i:(u,v)\in E, u<j_0, v\geq j_0\}.
\eeqn
Then by Lemma~\ref{lem:j0} and the Markov property, we get
\beqn
\nonumber
\h_{ij}(y,w,w') &=&\\
&&
\nonumber
\hspace{-2cm}
{\ts\oo2}
\sum_{x\bsub{C}}
\labs
\sum_{x\bsub{Z}}
\lparen{
\pr{X\bsub{C\cup Z}=x\bsub{C\cup Z}\gn X_i=w}
}
-    
\pr{X\bsub{C\cup Z}=x\bsub{C\cup Z}\gn X_i=w'}
\rparen\rabs
\\\label{eq:hZC}&&
\eeqn
(the sum indexed by $\{j_0,\ldots,n\}\setminus C$ marginalizes out).

Define $D=\{d_k:k=0,\ldots,|D|\}$ with 
$d_0=\depth_T(i)$,
$d_{|D|}=\depth_T(j_0)$
and $d_{k+1}=d_k+1$ for $0\leq k<|D|$.
For $d\in D$, 
let $I_d=T_i\cap \lev_{T}(d)$ 
and $G_d=(I_{d-1}+I_d,E_d)$ be the
bipartite graph consisting of the nodes in $I_{d-1}$ and $I_{d}$,
and the edges in $E$ joining them (note that 
$I_{d_0}=\{i\}$).

For $(u,v)\in E$, let $\vec A\supr{u,v}$ be the $|\X|\times|\X|$
matrix given by
\beq
\vec A\supr{u,v}_{x,x'} = p_{uv}(x\gn x')
\eeq
and note that $\nrm{\vec A\supr{u,v}}=\tha_{uv}$.
Then by the Markov property, for each 
$z\bsub{I_{d}}\in\X^{I_d}$ and
$x\bsub{I_{d-1}}\in\X^{I_{d-1}}$, 
$d\in D\setminus\{d_0\}$,
we have
\beq
\pr{X\rsub{I_{d}}=z\rsub{I_{d}}\gn X\rsub{I_{d-1}}=x\rsub{I_{d-1}}}
&=&
\vec A\supr{d}\bsub{z\rsub{I_{d}},x\rsub{I_{d-1}}},
\eeq
where
\beq
\vec A\supr{d} &=&
\TP_{(u,v)\in E_d} \vec A\supr{u,v}.
\eeq
Likewise, for 
$d\in D\setminus\{d_0\}$,
\beqn
\pr{X\rsub{I_{d}}=x\rsub{I_{d}}\gn X_i=w}
&=&
\sum_{x\rsub{I_1}'}\sum_{x\rsub{I_2}''}\cdots\sum_{x\rsub{I_{d-1}}\supr{d-1}}
\nonumber\\&&
\pr{X\rsub{I_{1}}=x\rsub{I_{1}}'\gn X_i=w}
\pr{X\rsub{I_{2}}=x\rsub{I_{2}}''\gn X\rsub{I_{1}}=x\rsub{I_{1}}'}
\cdots
\nonumber\\
&&\pr{X\rsub{I_{d}}=x\rsub{I_{d}}\gn X\rsub{I_{d-1}}=x\rsub{I_{d-1}}\supr{d-1}}
\nonumber\\
&=& 
\label{eq:AA}
(\vec A\supr d\vec A\supr{d-1}\cdots\vec A\supr{d_1})\bsub{x_{I_d},w}.
\eeqn
Define the 
(balanced)
$I_{d_1}$-tensor
\beqn
\label{eq:hvdef}
\vec h=
\vec A\supr{d_1}\bsub{\cdot,w} - \vec A\supr{d_1}\bsub{\cdot,w'},
\eeqn
the $I_{d_{|D|}}$-tensor
\beqn
\label{eq:fvdef}
 \vec f = \vec A\supr{d_{|D|}}\vec A\supr{d_{|D|-1}}\cdots\vec
A\supr{d_2}\vec h,
\eeqn
and $C_0,C_1,Z_0
\subset 
\{1,\ldots,n\}
$:
\beqn
C_0=C\cap I_{\depth_T(j_0)},
\qquad
C_1 = C\setminus C_0,
\qquad
Z_0 = I_{\depth_T(j_0)}\setminus C_0,
\eeqn
where $C$ and $Z$ are defined in (\ref{eq:ZC}).
For readability we will write $\P(x_U\gn\cdot)$ instead of
$\pr{X_U=x_U\gn\cdot}$ below; no ambiguity should arise.
Combining (\ref{eq:hZC}) and (\ref{eq:AA}), we have
\beqn
\h_{ij}(y,w,w') &=&
{\ts\oo2}
\sum_{x\rsub{C}}
\dsabs{
\sum_{x\rsub{Z}}
\paren{
\P(x\bsub{C\cup Z}\gn X_i=w)
-    
\P(x\bsub{C\cup Z}\gn X_i=w')
}}\\
&=&
{\ts\oo2}
\sum_{x\rsub{C_0}}
\sum_{x\rsub{C_1}}
\labs
\sum_{x\rsub{Z_0}}
\P(x\bsub{C_1}\gn x\bsub{Z_0})
\vec f\bsub{C_0\cup Z_0}
\rabs\\
&=&\nrm{\vec{Bf}}
\eeqn
where $\vec B$ is the
$|\X^{C_0\cup C_1}|\times|\X^{C_0\cup Z_0}|$ column-stochastic matrix
given by
$$ \vec B\bsub{x\rsub{C_0}\cup x\rsub{C_1},x'\rsub{C_0}\cup x_{Z_0}}
= \pred{x\rsub{C_0}=x'\rsub{C_0}} \P(x\rsub{C_1}\gn x\rsub{Z_0})$$
with the convention that $\P(x\rsub{C_1}\gn x\rsub{Z_0})=1$ if 
either of
$\{Z_0$,$C_1\}$
is empty. The claim now follows by reading off the results previously
obtained:
\beq
\begin{array}{rcllll}
\nrm{\vec{Bf}}
&\leq& \nrm{\vec{B}}\nrm{\vec{f}}
&&&\text{Eq. (\ref{eq:tv})}\\\\
&\leq&\nrm{\vec{f}}
&&&\text{Remark~\ref{rem:stochnorm}}\\\\
&\leq& \nrm{\vec h}\prod_{k=2}^{|D|} \nrm{\vec A\supr{d_k}}
&&&\text{Eqs. (\ref{eq:AB},\ref{eq:fvdef})}\\\\
&\leq& \prod_{k=1}^{|D|}
\inex\{\dsnrm{\vec A\supr{u,v}}:(u,v)\in E_{d_k}\}
&&&\text{Lemma~\ref{lem:TP}}.
\end{array}
\eeq
\enpf

\bepf[Proof of Theorem~\ref{thm:main}]
We will borrow 
the definitions from the proof of 
Lemma~\ref{lem:maintech}.
To upper-bound $\bar\h_{ij}$ we first bound
$\inex\{\dsnrm{\vec A\supr{u,v}}:(u,v)\in E_{d_k}\}$. Since 
$$ |E_{d_k}|\leq \width(T)\leq L$$
(because every node in $I_{d_k}$ has exactly one parent in
$I_{d_{k-1}}$) and
$$ \nrm{\vec A\supr{u,v}}=\tha_{uv} \leq \tha < 1,$$ we appeal to
Lemma~\ref{lem:inex} to obtain
\beqn
\inex\{\dsnrm{\vec A\supr{u,v}}:(u,v)\in E_{d_k}\}
&\leq&
1-(1-\tha)^L.
\eeqn
Now we must lower-bound the quantity
$h=\depth_T(j_0)-\depth_T(i)$. Since every level can have up to $L$ nodes,
we have
$$ j_0-i \leq hL $$
and so $h\geq\floor{(j_0-i)/L}\geq\floor{(j-i)/L}$.
\enpf

The calculations 
in Lemma~\ref{lem:maintech}
yield
considerably more information than the simple bound in (\ref{eq:mainbd}).
For example, suppose the tree $T$ has 
levels
$\{I_d:d=0,1,\ldots\}$
with the property that the levels are growing
at most linearly:
$$ |I_d|\leq cd $$
for some $c>0$. Let 
$d_i=\depth_T(i)$,
$d_j=\depth_T(j_0)$,
and $h=d_j-d_i$.
Then
\beq
j-i \leq j_0-i &\leq&
c\sum_{d_i+1}^{d_j} k \\
&=& \frac{c}{2}(d_j(d_j+1)-d_i(d_i+1))\\
&<& \frac{c}{2}((d_j+1)^2-d_i^2)\\
&<& \frac{c}{2}(d_i+h+1)^2
\eeq
so 
\hide{
There is no loss of generality in taking $i=1$ (since by
Lemma~\ref{lem:j0}) only the nodes in $T_i$ affect the value of
$\h_{ij}$. Then
$$ j-i \leq c\sum_{k=1}^h k = \frac{h(h+1)}{2} < c(h+1)^2/2, $$
so $$ h > \sqrt{2(j-i)/c}-1,$$
}
$$ h > \sqrt{2(j-i)/c}-d_i-1,$$
which yields the bound, via Lemma~\ref{lem:inex}\iif,
\beqn
\bar\h_{ij} &\leq& \prod_{k=1}^h 
\sum_{(u,v)\in E_k} \tha_{uv}.
\eeqn
Let
$\tha_k =\max\{\tha_{uv}: (u,v)\in E_k\}$;
then if
$ck\tha_k\leq\beta$ 
holds
for some 
$\beta\in\RR$,
this becomes
\beqn
\nonumber
\bar\h_{ij} &
\leq
& \prod_{k=1}^h (ck\tha_k)\\
\nonumber
&<&  \prod_{k=1}^{\sqrt{2(j-i)/c}-d_i-1} (ck\tha_k)\\
&\leq& 
\label{eq:thasqrt}
\beta^{\sqrt{2(j-i)/c}-d_i-1}.
\eeqn
This
is a non-trivial bound for trees with linearly growing
levels: recall that to bound $\nrm{\Delta}_\infty$
(\ref{eq:infnorm}), 
we must bound the series
$$ \sum_{j=i+1}^\infty \bar\h_{ij}.$$
By the limit comparison test with the series $\sum_{j=1}^\infty
1/j^2
$, we have that
$$ \sum_{j=i+1}^\infty \beta^{\sqrt{2(j-i)/c}-d_i-1}$$
converges for $\beta<1$.
Similar techniques may be applied when the level growth is
bounded by other slowly increasing functions.

\section{Discussion}
\label{sec:discussion}
We have presented a concentration of measure bound for
Markov tree processes; to our knowledge, this is the first such
result.\footnote{
In 
a 2003 paper,
Dembo et al.~\cite{dembo03} presented large deviation bounds for typed
Markov trees, which is a more general class of processes than the
Markov tree processes defined here. The techniques used and bounds
obtained in~\cite{dembo03} are of a rather different flavor than here;
this is not surprising since measure concentration and large
deviations, while pursuing similar goals, tend to use different
methods and state results that are often not immediately comparable.
}
In the simple case of the {\it contracting, bounded-width} Markov tree
processes (i.e., those for which $\width(T)\leq L<\infty$ and
$\sup_{u,v}\tha_{uv}\leq\tha<1$), the bound takes on a particularly
tractable form (\ref{eq:mainbd}), 
and
in the degenerate case 
$L=1$
it reduces to
the sharpest
known bound
for Markov chains. 
The techniques we develop extend well beyond the somewhat
restrictive contracting-bounded-width case, as demonstrated in the
calculation in (\ref{eq:thasqrt}).

The technical results in \S\ref{sec:proofmain}, particularly 
Lemma~\ref{lem:tvineq} and its generalizations, might be of
independent interest. It is hoped 
that these techniques will be extended to obtain concentration bounds
for larger classes of directed acyclic graphical models.

\section*{Acknowledgements}
I thank 
John Lafferty and
Kavita Ramanan for useful discussions and suggestions.


\end{document}